\documentclass{ifacconf}
\usepackage{subfigure}
\usepackage{amssymb, amsmath, amsfonts}
\usepackage{empheq}
\usepackage{caption, cite}
\usepackage{graphicx}      
\usepackage{natbib}        

\usepackage{algorithm}
\usepackage{algorithmic}

\usepackage{tikz}
\usetikzlibrary{shadows,arrows,positioning}

\pgfdeclarelayer{background}
\pgfdeclarelayer{foreground}
\pgfsetlayers{background,main,foreground}

\usetikzlibrary{external}
\usepackage{pgfplots}
\usepackage{cases}

 \newtheorem{theorem}{Theorem}
 \newtheorem{definition}{Definition}
 \newtheorem{lemma}{Lemma}

\newtheorem{remark}{Remark}
\newtheorem{assumption}{Assumption}

\newlength\figureheight
\newlength\figurewidth

\DeclareFontFamily{OT1}{pzc}{}
\DeclareFontShape{OT1}{pzc}{m}{it}{<-> s * [1.000] pzcmi7t}{}
\DeclareMathAlphabet{\mathpzc}{OT1}{pzc}{m}{it}

\newcommand{\R}{{\mathbb{R}}}
\newcommand{\E}{{\mathbb{E}}}

\newcommand{\x}{{\mathbf{x}}}
\newcommand{\bx}{{\mathbf{\bar{x}}}}
\newcommand{\sx}{{\mathbf{\sigma_x}}}
\newcommand{\y}{{\mathbf{y}}}
\newcommand{\sy}{{\mathbf{\sigma_{\y}}}}
\newcommand{\by}{{\mathbf{\bar{\y}}}}
\newcommand{\f}{{\mathbf{\nabla{f}}}}
\newcommand{\bbf}{{\mathbf{\nabla {\bar{f}}}}}
\newcommand{\lf}{{\nabla{f}}}
\newcommand{\ox}{{\Omega_{\x}}}
\newcommand{\obx}{{\Omega_{\bar{\x}}}}
\newcommand{\oy}{{\Omega_{\y}}}
\newcommand{\osx}{{\Omega_{\sigma_{\x}}}}
\newcommand{\osy}{{\Omega_{\sigma_{\y}}}}
\newcommand{\ex}{{\mathbf{\eta_{\x}}}}

\newcommand{\ey}{{\mathbf{\eta_{\y}}}}
\newcommand{\eby}{{\mathbf{\bar{\eta}_{\y}}}}
\newcommand{\lw}{{\hat{\lambda}}}
\newcommand{\llg}{{\bar{\lambda}_{G}}}
\newcommand{\llw}{{\rho_{w}}}
\newcommand{\lwi}{{\bar{\lambda}_{W-I}}}

\newcommand\addtag{\refstepcounter{equation}\tag{\theequation}}

\makeatletter

\newcommand{\Rmnum}[1]{\expandafter\@slowromancap\romannumeral #1@}
\makeatother

\begin{document}
\begin{frontmatter}

\title{\LARGE \bf{Compressed Differentially Private Distributed Optimization with Linear Convergence}} 


\author[First]{Antai Xie} 
\author[Second]{Xinlei Yi} 
\author[First]{Xiaofan Wang}
\author[Third]{Ming Cao}
\author[First]{Xiaoqiang Ren}

\address[First]{School of Mechatronic Engineering and Automation, Shanghai University, Shanghai, China. Emails: \{xatai,\,xfwang,\,xqren\}@shu.edu.cn.}
\address[Second]{Lab for Information $\&$ Decision Systems, Massachusetts Institute of Technology, Cambridge, MA 02139, USA. Email: xinleiyi@mit.edu.}
\address[Third]{Faculty of Science and Engineering, University of Groningen, Groningen, the Netherlands. Email: m.cao@rug.nl.}

\begin{abstract}                
 This paper addresses the problem of differentially private distributed optimization under limited communication, where each agent aims to keep their cost function private while minimizing the sum of all agents' cost functions. In response, we propose a novel Compressed differentially Private distributed Gradient Tracking algorithm (CPGT). We demonstrate that CPGT achieves linear convergence for smooth and strongly convex cost functions, even with a class of biased but contractive compressors, and achieves the same accuracy as the idealized communication algorithm. Additionally, we rigorously prove that CPGT ensures differential privacy. Simulations are provided to validate the effectiveness of the proposed algorithm.
\end{abstract}


\end{frontmatter}

\section{Introduction}
 In recent decades, the problem of distributed optimization in multi-agent systems has gained significant attention due to its wide applicability in various fields, such as sensor networks~large-scale machine learning~\cite[]{tsianos2012consensus},~\cite[]{dougherty2016extremum}, and online optimization~\cite[]{li2020distributed}. In a typical setup for distributed consensus optimization, the goal is to minimize the sum of all agents' local cost functions over a connected network, where each agent has access only to its own local cost function. This problem has been extensively studied in the literature, leading to the development of methods such as distributed (sub)gradient descent~\cite[]{nedic2009distributed}, EXTRA~\cite[]{shi2015extra}, and gradient tracking~\cite[]{qu2017harnessing}.
	 
The exchange of information over wireless networks may be vulnerable to attacks by malicious adversaries. Recent literature has shown that adversaries can obtain private training sets from shared gradients in machine learning~\cite[]{zhu2019deep}, potentially resulting in the exposure of sensitive data such as medical records and financial transactions. Thus, it is crucial to address privacy concerns in distributed optimization as a matter of urgency.
	 
To preserve the privacy of each agent,~\cite{huang2015differentially};~\cite{ zhu2018differentially};~\cite{ding2021differentially};~\cite{chen2021differential} proposed several differentially private distributed optimization algorithms by introducing the notion of differential privacy~\cite[]{dwork2008differential}, which has mathematically provable security properties. More specifically, \cite{huang2015differentially} proposed differentially private gradient descent by masking states with Laplacian noises. \cite{zhu2018differentially} extended the above results to time-varying directed networks. However, due to the limitation of decaying stepsizes, both approaches can only guarantee sublinear convergence. To this end, \cite{ding2021differentially} achieved both linear convergence and differential privacy by simultaneously adding noise to states and directions and using constant stepsize. \cite{chen2021differential} further considered the case of directed graphs. Furthermore, another main approach to reach privacy-preserving is encryption. For example, \cite{lu2018privacy} proposed a privacy-preserving distributed optimization method using homomorphic encryption. Although these encryption-based methods can enable the solutions to converge to the exact optimal, they require a considerable amount of computing resources.
	
Most of the aforementioned approaches investigated the privacy-preserving distributed optimization algorithms under the idealized communication network. In practice, due to the limited communication bandwidth, it is necessary to consider compressed information. To this end, various research results have been proposed. \cite{alistarh2017qsgd} proposed communication-efficiency stochastic gradient descent algorithms by using an unbiased compressor. \cite{kajiyama2020linear} and \cite{liao2022compressed} achieved linear convergence by combining the gradient tracking algorithm and a compressor with bounded absolute compression error and a class of compressor with bounded relative compression error, respectively. \cite{xiong2021quantized} extended the approach in~\cite{kajiyama2020linear} to directed graphs.

 Compressed information offers many advantages, such as reducing communication costs, making it natural to consider combining communication compression with privacy preservation. To this end, \cite{wang2022quantization} proposed a compressed differentially private gradient descent algorithm that incorporates both compression and privacy preservation. However, they did not provide the linear convergence analysis.
	 
In this paper, inspired by~\cite{ding2021differentially}, we propose a compressed differentially private distributed gradient tracking algorithm, which achieves linear convergence and preserves differential privacy. The main contributions of this work are summarized as follows:
\begin{enumerate}
\item For a class of biased but contractive compressors, we propose a novel Compressed differentially private Gradient Tracking algorithm (CPGT). We show that CPGT achieves linear convergence (Theorem~\ref{theo:convergence}) and has the same accuracy as that of the algorithm over idealized communication network~\cite[]{ding2021differentially}.
\item We show that CPGT preserves differential privacy for the local cost function of each agent (Theorem~\ref{theo:privacy}). It is worth noting that CPGT is effective for a class of compressors. Different from~\cite{wang2022quantization}, CPGT does not rely on any specific compressor.
\end{enumerate}

The remainder of this paper is organized as follows. In Section~\ref{sec:Problemsetup}, we introduce the necessary definitions and formulate the considered problem. The CPGT is proposed in Section~\ref{sec:Algorithm}, and the convergence and privacy of the proposed CPGT are analyzed. Section~\ref{simulation} provides a numerical example to illustrate the results. Finally, Section~\ref{conclusion} concludes the paper. 
	
\emph{Notations}: $\R$ ($\mathbb{R}_{+}$) is the set of (positive) real numbers. $\mathbb{Z}$ is the set of integers and $\mathbb{N}$ the set of nature numbers. $\mathbb{R}^n$ and $\mathbb{R}^{n\times d}$ are the set of $n$ dimensional vectors and $n\times d$ dimensional matrices with real values, respectively. The transpose of a matrix $P$ is denoted by $P^T$, and we use $[P]_{ij}$ to denote the element in $i$-th row and $j$-th column. The all-ones and all-zeros column vector are denoted by $\mathbf{1}$ and $\mathbf{0}$, respectively. The identity matrix is denoted by $I$. We then introduce a stacked matrix: for a matrix $\x\in\R^{n\times d}$, $\mathbf{\bar{\x}}\triangleq\frac{\mathbf{1}\mathbf{1}^T}{n}\mathbf{\x}$. $\vert\cdot\vert,\Vert\cdot\Vert,\Vert\cdot\Vert_F$ denote the absolute value, $l_2$ norm, and Frobenius norm, respectively. For a matrix $W$, we use $\bar{\lambda}_W$ to denote its spectral radius. For a given constant $\theta>0$, $\text{Lap}(\theta)$ is the Laplace distribution with the probability function $f_L(x,\theta)=\frac{1}{2\theta}e^{\frac{\vert x\vert}{\theta}}$. For any vector $\eta=[\eta_1,\dots,\eta_d]^T\in\R^d$, we say that $\eta\sim\text{Lap}_d(\theta)$ if each component $\eta_i\sim\text{Lap}(\theta)$, $i=1,\dots,d$. Furthermore, we use $\E(x)$ and $\mathbb{P}(x)$ to denote the expectation and probability of a random variable $x$, respectively. In addition, if $x\sim\text{Lap}(\theta)$, we have $\E(x^2)=\theta^2$ and $\E(\vert x\vert)=\theta$.

\section{Preliminaries and Problem Formulation} \label{sec:Problemsetup}

\subsection{Distributed Optimization}

We consider a network of $n$ agents, where each agent has a private convex cost function $f_i:\mathbb{R}^d\mapsto\mathbb{R}$. All agents solve the following optimization problem cooperatively:
\begin{align}\label{P1}
	\min_{x\in\mathbb{R}^d}\sum_{i=1}^n f_i(x),
\end{align}
where $x$ is the global decision variable, which is not known by each agent and can only be estimated locally. More precisely, we assume each agent $i$ maintains a local estimate $x_i(k)\in\mathbb{R}^d$ of $x$ at time step $k$ and use $\lf_i(x_i(k))$ to denote the gradient of $f_i$ with respect to $x_i(k)$. Moreover, we make the following assumptions on the local cost functions~$f_i$:
\begin{assumption}\label{as:smoothandsconvex}
	Each local cost function $f_i$ is $\mu$-strongly convex and $L$-smooth, where $0<\mu\leq L$. That is, for any~$x,y\in\mathbb{R}^d$,
	\begin{align*}
		&f_i(y)\geq f_i(x)+\lf_i(x)^T(y-x)+\frac{\mu}{2}\Vert x-y\Vert^2,\\
		&\Vert \lf_i(x)-\lf_i(y)\Vert\leq L\Vert x-y\Vert.
	\end{align*}
\end{assumption}

Assumption~\ref{as:smoothandsconvex} is standard for linear convergence in distributed optimization, e.g.,~\cite{xu2017convergence,qu2017harnessing,kajiyama2020linear,xiong2021quantized}. Furthermore, under Assumption~\ref{as:smoothandsconvex}, Problem~\eqref{P1} has a unique optimal solution $x^*$~\cite[]{boyd2004convex}.

\subsection{Basics of Graph Theory}
	
The exchange of information between agents is captured by an undirected graph $\mathcal{G}(\mathcal{V}, \mathcal{E})$ with $n$ agents, where $\mathcal{V}=\{1, 2, \ldots, n\}$ is the set of the agents' indices and $\mathcal{E} \subseteq \mathcal{V} \times \mathcal{V}$ is the set of edges. The edge $(i,j)\in\mathcal{E}$ if and only if agents~$i$ and~$j$ can communicate with each other. Let $W=[w_{ij}]_{n\times n}\in\mathbb{R}^{n\times n}$ be the adjacency matrix of $\mathcal{G}$, namely $w_{ij}>0$ if $(i,j)\in\mathcal{E}$ or $i=j$, and $w_{ij}=0$ otherwise. We use $\mathcal{N}_i=\{j\in\mathcal{V}|~(i,j)\in\mathcal{E}\}$ to denote the neighbor set of agent $i$.
	
\begin{assumption}\label{as:strongconnected}
	The undirected graph $\mathcal{G}(\mathcal{V}, \mathcal{E})$ is connected and $W$ is a doubly stochastic matrix, i.e., $\mathbf{1}^TW=\mathbf{1}^T$ and $W\mathbf{1}=\mathbf{1}$.
\end{assumption}

\subsection{Compression Method}
Due to the limited communication channel capacity, we consider the situation where agents compress the information before sending it. More specifically, for any $x\in\R^d$, we consider a class of stochastic compressors $C(x,\varrho)$ and use $f_c(x,\varrho)$ to denote the corresponding probability density functions, where $\varrho$ is a random perturbation variable. $C(x,\varrho)$ can be simplified to $C(x)$ when the distribution of $\varrho$ is given. We then introduce the following assumption.

\begin{assumption}\label{as:compressor}
For some $\varphi\in[0,1)$, the compressor $C(\cdot):\mathbb{R}^d\mapsto\mathbb{R}^d$ satisfies 
\begin{align*}
    \mathbb{E}_C\begin{bmatrix}
        \Vert C(x)-x\Vert^2
    \end{bmatrix}\leq\varphi\Vert x\Vert^2, \forall x\in\mathbb{R}^d.\addtag\label{eq:propertyofcompressors}
\end{align*}
\end{assumption}
From~\eqref{eq:propertyofcompressors} and the Jensen's inequality, one obtains that
\begin{align}\label{eq:propertyofcompressors1}
    \mathbb{E}_C\begin{bmatrix}
        \Vert C(x)-x\Vert
    \end{bmatrix}\leq\sqrt{\varphi}\Vert x\Vert, \forall x\in\mathbb{R}^d.
\end{align}
\begin{remark}
Compressors under Assumption~\ref{as:compressor} are quite common, see e.g. \cite{reisizadeh2019exact,koloskova2019decentralized,taheri2020quantized}. Note that, unlike~\cite{wang2022quantization}, Assumption~\ref{as:compressor} does not require compressors to be unbiased.
\end{remark}

\subsection{Differential Privacy}

To evaluate the privacy performance, we adopt the notion of $(\epsilon,\delta)$-differential privacy for the distributed optimization, which has recently been studied in~\cite{huang2015differentially};~\cite{ding2021differentially}. Specifically, we introduce the following definitions.

\begin{definition}\label{def:adjacent}
\textbf{($\delta$-adjacent~\cite[]{ding2021differentially})} Given $\delta>0$, two function sets $\mathcal{S}^{(1)}=\{f_i^{(1)}\}_{i=1}^n$ and $\mathcal{S}^{(2)}=\{f_i^{(2)}\}_{i=1}^n$ are said to be $\delta$-adjacent if there exists some $i_0\in\{1,2,\dots,n\}$ such that
\begin{align*}
    f_i^{(1)}=f_i^{(2)}~\forall i\neq i_0,~\text{and}~D(f_{i_0}^{(1)},f_{i_0}^{(2)})\leq\delta,
\end{align*}
where $D(f_{i_0}^{(1)},f_{i_0}^{(2)})\triangleq\max_{x\in\mathbb{R}^d}\Vert \lf_{i_0}^{(1)}(x)-\lf_{i_0}^{(2)}(x)\Vert$ represents the distance between two functions $f_{i_0}^{(1)}$ and $f_{i_0}^{(2)}$.
\end{definition}

\begin{definition}\label{def:differentialprivacy}
\!\!\textbf{(Differential privacy~\cite[]{chen2021differential}\!)} Given $\delta,\epsilon>0$ and a randomized mechanism $M$, for any two $\delta$-adjacent function sets $\mathcal{S}^{(1)}$ and $\mathcal{S}^{(2)}$, and any observation $\mathcal{H}\subseteq\text{Range}(M)$, the randomized mechanism $M$ keeps $\epsilon$-differential privacy if
\begin{align}
\mathbb{P}\{M(\mathcal{S}^{(1)})\in\mathcal{H}\}\leq e^\epsilon\mathbb{P}\{M(\mathcal{S}^{(2)})\in\mathcal{H}\}.
\end{align}
\end{definition}

Definition~\ref{def:differentialprivacy} shows that the randomized mechanism $M$ is differentially private if for any pair of $\delta$-adjacent function sets, the probability density functions of their observations are similar. Intuitively, it is difficult for an adversary to distinguish between two $\delta$-adjacent function sets merely by observations if the corresponding mechanism $M$ is differential private.

\section{Main Results}\label{sec:Algorithm}
In this section, we provide the Compressed differentially private Gradient Tracking algorithm (CPGT), which is shown in Algorithm~\ref{Al:CPGT}.
\begin{algorithm}[]
	\caption{CPGT Algorithm }
	\label{Al:CPGT}
	\begin{algorithmic}[1]
		\STATE \textbf{Input:} Stopping time $K$, adjacency matrix $W$, and positive parameters $\alpha$, $\gamma$, $d_{\eta_{x_i}}$, $d_{\eta_{y_i}}$, $q_i$, $\forall i\in\mathcal{V}$.
		\STATE \textbf{Initialization:} Each ~$i\in\mathcal{V}$ chooses arbitrarily $x_i(0)\in\mathbb{R}^d$, $x^c_i(-1)=\bf{0}$, $y^c_i(-1)=\bf{0}$, and computes $y_i(0)=\lf_i(x_i(0))$.
		\FOR{$k=0,1,\dots,K-1$}
		\FOR {for $i\in\mathcal{V}$ in parallel} 
		\STATE Generate Laplace noises $\eta_{x_i}(k)\sim\text{Lap}_d(d_{\eta_{x_i}}q_i^k)$ and $\eta_{y_i}(k)\sim\text{Lap}_d(d_{\eta_{y_i}}q_i^k)$.
		\STATE Obtain $x_i^a(k)$ and $y_i^a(k)$ from~\eqref{eq:addnoisex} and~\eqref{eq:addnoisey}, respectively.
		\STATE Compute $C(x_i^a(k)-x_i^c(k-1))$ and $C(y_i^a(k)-y_i^c(k-1))$, then broadcast them to its neighbors $\mathcal{N}_i$.
		\STATE Receive $C(x_j^a(k)-x_j^c(k-1))$, and $C(y_j^a(k)-y_j^c(k-1))$ from $j\in\mathcal{N}_i$.
		\STATE Update $x_j^c(k)$ and $y_j^c(k)$ from~\eqref{citerationx} and~\eqref{citerationy}, respectively.
		\STATE Update $x_i(k+1)$ and $y_i(k+1)$ from~\eqref{iterationx} and~\eqref{iterationy}, respectively.
		\ENDFOR
		\ENDFOR
		\STATE \textbf{Output:} \{$x_i(K)$\}.
	\end{algorithmic}
\end{algorithm}
\subsection{Algorithm Description}
The proposed CPGT is inspired by DiaDSP~\cite[]{ding2021differentially}. We assume each agent $i\in\mathcal{V}$ maintains an estimate $x_i(k)$ and an auxiliary variable $y_i(k)$ for tracking the global gradient. To enable differential privacy, each agent $i$ broadcasts the noise added information $x_i^a(k)$ and $y_i^a(k)$ to its neighbors $\mathcal{N}_i$ per step, where
\begin{align*}
&x_i^a(k)=x_i(k)+\eta_{x_i}(k),\addtag\label{eq:addnoisex}\\
&y_i^a(k)=y_i(k)+\eta_{y_i}(k),\addtag\label{eq:addnoisey}
\end{align*}
with $\eta_{x_i}(k)$ and $\eta_{y_i}(k)$ are Laplace noises. Similar to the DiaDSP Algorithm~\cite[]{ding2021differentially}, we set $\eta_{x_i}\sim\text{Lap}_d(d_{\eta_{x_i}}q_i^k)$ and $\eta_{y_i}\sim\text{Lap}_d(d_{\eta_{y_i}}q_i^k)$, $\forall i\in\mathcal{V}$, where $d_{\eta_{x_i}}>0,~d_{\eta_{y_i}}>0$, and $0<q_i<1$. After the information exchange, agent $i$ performs the following updates:
\begin{align}
	&x_i(k+1)=\sum_{j=1}^n w_{ij}x_j^a(k)-\alpha y_i(k),\\
	&y_i(k+1)=\sum_{j=1}^n w_{ij}y_j^a(k)+\lf_i(x_i(k+1))-\lf_i(x_i(k)),
\end{align}
where the stepsize $\alpha$ is a constant and the initial value $y_i(0)=\lf_i(x_i(0)),~\forall i\in\mathcal{V}$. To improve the communication efficiency, we use the compressed information $x_i^c(k)$, $y_i^c(k)$ to replace $x_i^a(k)$ and $y_i^a(k)$, respectively. Then, we design the updates of agent $i\in\mathcal{V}$ as follows:
\begin{align*}
	&~x_i(k+1)=x_i^a(k)+\gamma\sum_{j=1}^n w_{ij}(x_j^c(k)-x_i^c(k))-\alpha y_i(k),\addtag\label{iterationx}\\
	&~y_i(k+1)=y_i^a(k)+\gamma\sum_{j=1}^n w_{ij}(y_j^c(k)-y_i^c(k))\\
	&~~~~~~~~~~~~~~~+\lf_i(x_i(k+1))-\lf_i(x_i(k)),\addtag\label{iterationy}
\end{align*}
where
\begin{align*}
	&x_i^c(k)=x_i^c(k-1)+C(x_i^a(k)-x_i^c(k-1)),\addtag\label{citerationx}\\
	&y_i^c(k)=y_i^c(k-1)+C(y_i^a(k)-y_i^c(k-1)),\addtag\label{citerationy}
\end{align*}
with $\gamma$ being a positive parameter. We assume that $x_i^c(k)=\mathbf{0}$ and $y_i^c(k)=\mathbf{0}$ for $k<0,~\forall i\in\mathcal{V}$. Let $W_\gamma\triangleq(1-\gamma)I+\gamma W$, \eqref{iterationx} and \eqref{iterationy} can be rewritten into the following matrix form
\begin{align*}
	&\mathbf{x}(k+1)=W_\gamma(\x(k)+\ex(k))+\gamma(W-I)\sx(k)-\alpha\y(k),\addtag\label{iterationx1}\\
	&\mathbf{y}(k+1)=W_\gamma(\y(k)+\ey(k))+\gamma(W-I)\sy(k)\\
	&~~~~~~~~~~~~+\f(\x(k+1))-\f(\x(k)),\addtag\label{iterationy1}
\end{align*}
where $\x(k)\triangleq[x_1(k),x_2(k),\dots,x_n(k)]^T\in\R^{n\times d}$, $\y(k)\triangleq[y_1(k),x_2(k),\dots,y_n(k)]^T\in\R^{n\times d}$, $\sx(k)\triangleq[x_1^c(k)-x_1^a(k),$ $\dots,x_n^c(k)-x_n^a(k)]^T\in\R^{n\times d}$, $\sy(k)\triangleq[y_1^c(k)-y_1^a(k),\dots,$ $y_n^c(k)-y_n^a(k)]^T\in\R^{n\times d}$, $\f(\mathbf{x}(k))\triangleq[\lf_1(x_1(k)),\dots,$ $\lf_n(x_n(k))]^T\in\R^{n\times d}$, $\ex(k)\triangleq\left[\eta_{x_1}(k),\dots,\eta_{x_n}(k)\right]^T\in\R^{n\times d}$ and $\ey(k)\triangleq\left[\eta_{y_1}(k),\dots,\eta_{y_n}(k)\right]^T\in\R^{n\times d}$.

\subsection{Convergence Analysis of CPGT}
In this section, we first show linear convergence of CPGT under the compressors satisfying Assumption~\ref{as:compressor}. Second, the differential privacy of all cost functions is proved under the CPGT Algorithm. We would like to point out that CPGT has the same convergence accuracy as the algorithm with idealized communication, i.e., DiaDSP~\cite[]{ding2021differentially}.

Let $\Theta(k)\triangleq[\ox(k), \obx(k), \oy(k), \osx(k), \osy(k)]^T$, where $\ox(k)=\Vert \x(k)-\bx(k)\Vert_F$, $\obx(k)=\Vert \bx(k)-\x^\infty\Vert_F$, $\oy(k)=\Vert \y(k)-\by(k)\Vert_F$, $\osx(k)=\Vert \sx(k)\Vert_F$, and $\osy(k)=\Vert \sy(k)\Vert_F$, with $\x^\infty\in\R^{n\times d}$ being given in Lemma~\ref{lemma:linearinequalities}. The following lemma constructs a linear system of inequalities that is related to $\Theta(k)$.

\begin{lemma}\label{lemma:linearinequalities}
	Let $\bar{q}=\max_i\{q_i\}$ and $\bar{d_\eta}=\max_i\{d_{\eta_{x_i}},d_{\eta_{y_i}}\}$. Suppose Assumptions~\ref{as:smoothandsconvex}--\ref{as:compressor} hold. Under Algorithm~\ref{Al:CPGT}, if $\alpha<\frac{1}{\mu+L}$, we have the following linear inequalities:
	\begin{align}\label{eq:linearinequalities}
		\mathbb{E}[\Theta(k+1)]\preceq G\mathbb{E}[\Theta(k)]+\vartheta \bar{q}^k\bar{d_\eta},
	\end{align}
where the notation $\preceq$ means element-wise less than or equal to, the elements of matrix $G\in\mathbb{R}^{5\times 5}$ and constant vector $\vartheta\in\R^{5\times 1}$ are given in Appendix~\ref{app-linearinequalities}, and $\x^\infty=\mathbf{1}(x^\infty)^T$ with $x^\infty\in\R^d$ satisfying
\begin{align}\label{eq:xinfty}
	\sum_{i=1}^n\lf_i(x^\infty)=-\sum_{i=1}^n\sum_{k=0}^\infty \eta_{y_i}(k).
\end{align}
\end{lemma}
\begin{pf}
	See Appendix~\ref{app-linearinequalities}.
\end{pf}
\vspace{-0.2em}
In light of Lemma~\ref{lemma:linearinequalities}, we know that CPGT can linearly converge to $\x^\infty$ if the spectral radius of matrix $G$ is strictly less than $1$, i.e., $\llg<1$. Hence, we establish linear convergence of CPGT by the following lemma.
\vspace{-0.4em}
\begin{lemma}\label{lemma:convergence}
Suppose Assumptions~\ref{as:smoothandsconvex}--\ref{as:compressor} hold. Under Algorithm~\ref{Al:CPGT}, $\llg\leq 1-\frac{m}{2}\alpha\mu$ holds for some $m\in(0,1)$, if parameter $\gamma$ and stepsize $\alpha$ satisfy
\begin{align*}
	&\gamma\leq\min\bigg\{\frac{(1-\sqrt{\varphi})\zeta_5}{2\sqrt{\varphi}\lwi(2\zeta_1+\zeta_2+2\zeta_3+\zeta_4+\zeta_5)},\\	&~~~~~~~~~~~~~\frac{(1-\sqrt{\varphi})L}{\lwi m\mu},1\bigg\},\addtag\label{eq:boundofgamma}\\
	&\alpha\leq\min\bigg\{\frac{\gamma(1-\llw)\zeta_1}{L(2\zeta_3+m\zeta_1)},\frac{\gamma(1-\llw)\zeta_3}{L(2\zeta_1+2\zeta_2+(2+m)\zeta_3)},\\
	&~~~~~~~~~~~~~\frac{\gamma\lwi}{2L}\bigg\},\addtag\label{eq:boundofstepsize}\\
	&\zeta_1+\zeta_4+\zeta_5\leq\frac{(1-\llw)\zeta_3}{2\lwi},\zeta_5\leq\zeta_4,\zeta_4\leq\frac{(1-\llw)\zeta_1}{2\lwi},\\
	&\zeta_1\leq\frac{m\mu\zeta_2}{2L},\addtag\label{eq:relationbetwzeta}
\end{align*}
where $\zeta_1,\dots,\zeta_5$ are some positive constants and~$\llw=\bar{\lambda}_{W-\frac{\mathbf{1}\mathbf{1}^T}{n}}$.
\end{lemma}
\vspace{-0.5em}
\begin{pf}
	See Appendix~\ref{app-convergence}.
\end{pf}
\vspace{-0.2em}
The following theorem shows that CPGT can linearly converge to the solution $\x^\infty$ in the mean by taking some concrete values for $\zeta_1,\dots,\zeta_5$.
\begin{theorem}\label{theo:convergence}
Suppose Assumptions~\ref{as:smoothandsconvex}--\ref{as:compressor} hold. Under Algorithm~\ref{Al:CPGT}, for some $m\in(0,1)$, if parameter $\gamma$ and stepsize $\alpha$ satisfy
\begin{align*}
	&\gamma\leq\min\left\{\frac{\kappa(1-\llw)^2 m\mu}{4\sqrt{\varphi}s_1},\frac{\kappa L}{m\mu},1\right\},\addtag\label{eq:boundofgamma1}\\
	&\alpha\leq\min\left\{\frac{m\mu(1-\llw)^2}{s_2},\frac{\lwi(1-\bar{q})}{2}\right\}\frac{\gamma}{L},\addtag\label{eq:boundofstepsize1}
\end{align*}
where $\kappa=\frac{1-\sqrt{\varphi}}{\lwi}$, $s_1=4\lwi^2m\mu+(1-\llw)^2m\mu+(6m\mu+2L)\lwi(1-\llw)$, and $s_2=(6m\mu+2m^2\mu+4L)(1-\llw)+2m\mu\lwi(2+m)$, then, we have $\E[\Vert\x(k)-\x^\infty\Vert_F]=\mathcal{O}((1-\frac{m}{2}\alpha\mu)^k)$.
\end{theorem}
\begin{pf}
Let $\zeta_1=1,~\zeta_2=\frac{2L}{m\mu},~\zeta_3=\frac{2\lwi-2\llw+2}{1-\llw},~\zeta_4=\zeta_5=\frac{1-\llw}{2\lwi}$, it is easy to verify that \eqref{eq:relationbetwzeta} holds. Then we derive \eqref{eq:boundofgamma}--\eqref{eq:boundofstepsize} from \eqref{eq:boundofgamma1}--\eqref{eq:boundofstepsize1} by substituting $\zeta_1,\dots,\zeta_5$ with the aforementioned values. From Lemmas~\ref{lemma:linearinequalities} and~\ref{lemma:convergence}, let $h=1-\frac{m}{2}\alpha\mu$, we have
\vspace{-0.2em}
\begin{align*}
	\mathbb{E}[\Vert\Theta(k+1)\Vert]&\leq h^k\mathbb{E}[\Vert\Theta(0)\Vert]+\vartheta\bar{d_\eta}\sum_{t=0}^{k-1}h^{k-1-t}\bar{q}^t\\
	&\leq h^k\mathbb{E}[\Vert\Theta(0)\Vert]+\vartheta\bar{d_\eta}h^{k-1}\sum_{t=0}^{k-1}\left(\frac{\bar{q}}{h}\right)^t\\
	&\leq (\mathbb{E}[\Vert\Theta(0)\Vert]+\frac{\vartheta\bar{d_\eta}}{h-\bar{q}})h^k.
\end{align*}
Since $\E[\Vert\x(k)-\x^\infty\Vert_F]\leq\E[\Vert\x(k)-\bx(k)\Vert_F]+\E[\Vert\bx(k)-\x^\infty\Vert_F]$, one obtains that
\begin{align}
	\E[\Vert\x(k)-\x^\infty\Vert_F]=\mathcal{O}(h^k),
\end{align}
which completes the proof.
\end{pf}
\begin{remark}
Noting that it is impossible to achieve both differential privacy and exact convergence simultaneously~(see \cite[Proposition 1]{ding2021differentially}). As shown in~\eqref{eq:xinfty}, the distance between convergence point $x^\infty$ and the optimal solution $x^*$ is affected by the sum of noise added to the gradients. Furthermore, convergence point $x^\infty$ of CPGT is the same as DiaDSP~\cite[]{ding2021differentially}, which means that the proposed CPGT has the same accuracy as the algorithm with idealized communication.
\end{remark}

\subsection{Differential Privacy}
In this section, we show that the differential privacy of all cost functions can be preserved under CPGT. 

We use $\mathcal{H}_k$ to denote the information transmitted between agents at time step $k$, i.e., $\mathcal{H}_k=\{C(x_i^a(k)-x_i^c(k-1)),C(y_i^a(k)-y_i^c(k-1))|~\forall i\in\mathcal{V}\}$. Without loss of generality, we assume the adversary aims to infer the cost function of agent $i_0$. Consider any two 
$\delta$-adjacent function sets $\mathcal{S}^{(1)}$ and $\mathcal{S}^{(2)}$, and only the cost function $f_{i_0}$ is different between the two sets, i.e., $f_{i_0}^{(1)}\neq f_{i_0}^{(2)}$ and $f_{i}^{(1)}= f_{i}^{(2)},~\forall i\neq i_0$.
\begin{assumption}\label{as:distance}
    For any $x_1,~x_2\in\R^d$, we have
    \begin{align*}
    \lf_{i_0}^{(1)}(x_1)-\lf_{i_0}^{(1)}(x_2)=\lf_{i_0}^{(2)}(x_1)-\lf_{i_0}^{(2)}(x_2).
    \end{align*}
\end{assumption}
\begin{remark}
    Assumption~\ref{as:distance} is necessary for granting the differential privacy in infinite time, and the same assumption can be founded in the literature~\cite[]{ding2021differentially}.
\end{remark}
\begin{theorem}\label{theo:privacy}
	Suppose Assumptions~\ref{as:smoothandsconvex}--\ref{as:distance} hold. CPGT preserves the $\epsilon_{i_0}$-differential privacy for any agent $i_0$'s cost function if $\alpha<\frac{1}{2L}$ and $q_{i_0}\in(\frac{\alpha L+\sqrt{\alpha^2L^2+4\alpha L}}{2},1)$, where $\epsilon_{i_0}$ is given by
\begin{align}
	\epsilon_{i_0}=\frac{\tau_{i_0} q_{i_0}^2\delta}{q_{i_0}^2-\alpha L-q_{i_0}\alpha L},~~~~\forall i_0\in\mathcal{V},
\end{align}
with $\tau_{i_0}=\frac{\alpha}{d_{\eta_{x_{i_0}}}}+\frac{1}{d_{\eta_{y_{i_0}}}}$.
\end{theorem}
\begin{pf}
	See Appendix~\ref{app-privacy}.
\end{pf}
\begin{remark}
Different from the approach in~\cite{wang2022quantization}, CPGT is effective for a class of compressors and achieves linear convergence.
\end{remark}

\section{simulation}\label{simulation}
In this section, simulations are given to verify the validity of CPGT. We consider a distributed estimation problem with $n=6$ agents and they communicate on a connected undirected graph, whose topology is shown in Fig.~\ref{Fig:graph}. Specifically, we assume each agent $i$ aims to estimate the unknown parameter $x\in\R^{10}$ and has noisy measurements $b_i=A_ix+e_i$, where $A_i\in\R^{6\times 10}$ is a nonsingular measurement matrix of agent $i$ and $e_i\in\R^{6}$ is a noise vector. This problem can be reformulated as the following form
\begin{align*}
	\min_{x}f(x)=\frac{1}{6}\sum_{i=1}^6\Vert A_ix-b_i \Vert^2.
\end{align*}
In this example, the elements of matrix $A_i$ and vector $b_i$ are randomly generated by Gaussian distributed $N(0,1)$. In addition, the initial value of each agent $x_i(0)$ is randomly chosen in $[0,1]^{10}$.
\begin{center}
	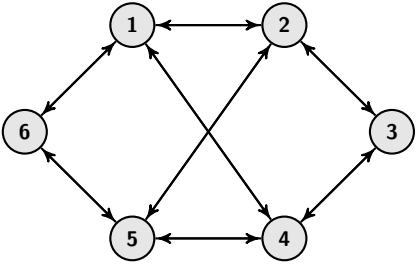
\begin{figure}
		\center
		\usetikzlibrary{positioning, fit, calc}

\begin{tikzpicture}[->,>=stealth',shorten >=1pt,auto,node distance=2cm,
thick,main node/.style={circle,fill=gray!20,draw,font=\sffamily\small\bfseries}]

\node[main node] (1) {1};
\node[main node] (2) [ right of=1] {2};
\node[main node] (3) [ below right of=2] {3};
\node[main node] (4) [ below left of=3] {4};
\node[main node] (5) [ left of=4] {5};
\node[main node] (6) [ below left of=1] {6};

%
%

\path[every node/.style={font=\sffamily\small}]
(1) edge  node[] {} (4)
edge  node[left] {} (2)	  
edge   node[above]{} (6)
(2) edge  node[left] {} (5)
edge  node[above]{} (1)
edge  node[left] {} (3)
(3) edge  node[left] {} (4)
edge  node[left] {} (2)
(4) edge  node[above]{} (3)
edge  node[above]{} (1)	
edge  node[left] {} (5)
(5) edge  node[left] {} (2)
edge node[left] {} (4)
edge node[left] {} (6)
(6) edge node[above]{} (5)
edge node[above]{} (1);
\end{tikzpicture}
		\caption{A connected undirected graph consisting of 6 agents.}
		\label{Fig:graph}
	\end{figure}
\end{center}
We consider the following two biased but contractive compressors:
\begin{itemize}
	\item Greedy (Top-$k$) quantizer~\cite[]{beznosikov2020biased}:
	\begin{align*}
		C_1(x):=\sum_{i_s=1}^kx_{(i_s)}e_{i_s},
	\end{align*}
where $x_{(i_s)}$ is the $i_s$-th coordinate of $x$ with $i_1,\dots,i_k$ being the indices of the largest $k$ coordinates in magnitude of $x$, and $e_1,\dots,e_d$ are the standard unit basis vectors in $\R^d$.
    \item biased $b$-bits quantizer~\cite[]{koloskova2019decentralized}:
    \begin{align*}
    	C_2(x):=\frac{\Vert x\Vert}{\xi}\cdot \text{sign}(x)\cdot 2^{-(b-1)}\circ \left\lfloor\frac{2^{(b-1)\vert x\vert}}{\Vert x\Vert} +u\right\rfloor,
    \end{align*}
where $\xi=1+\min\{\frac{d}{2^{2(b-1)}},\frac{\sqrt{d}}{2^{(b-1)}}\}$, $u$ is a random dithering vector uniformly sampled from $[0,1]^d$, $\circ$ is the Hadamard product, and $\text{sign}(\cdot)$, $|\cdot|$, $\lfloor\cdot\rfloor$ are the element-wise sign, absolute and floor functions, respectively. 
\end{itemize}	
As pointed out in~\cite{koloskova2019decentralized}, both of the above compressors satisfy Assumption~\ref{as:compressor}. Specifically, we choose $k=2$ and $b=2$ in the following simulations. 

\begin{figure}
	\centering
	\includegraphics[width=1\linewidth]{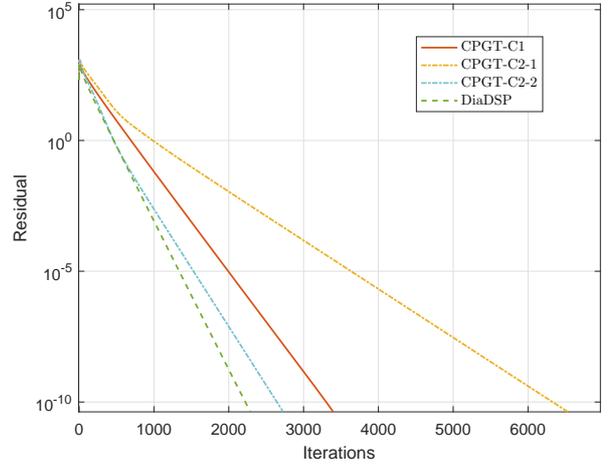}
	\caption{The evolution of residual under CPGT and uncompressed method~DiaDSP.}
	\label{fig:convergence}
\end{figure}

We first verify the convergence rate of CPGT with different compressors. We set $d_{\eta_{x_i}}=d_{\eta_{x_i}}=d_{\eta}$ and $q_i=q$, $\forall i\in\mathcal{V}$ and parameters of different algorithms are given in Table~\ref{tab:parameter}. The residual is computed, defined by $\text{residual}\triangleq\frac{1}{T}\sum_{t=1}^{T}\Vert\x_t(k)-\x^\infty\Vert_F$ with $T=1000$ and $\x_t(k)$ being the state $\x$ at time step $k$. Fig.~\ref{fig:convergence} shows that $\x(k)$ linearly converges to the point $\x^\infty$ under CPGT with different constant stepsizes and compressors. Moreover, the convergence rate of CPGT can be close to that of DiaDSP~\cite[]{ding2021differentially} with suitable parameters and compressors.
    \begin{table}[]
	\centering
	\begin{tabular}{lc c c c c}\hline
		Algorithm & Compressor & $\gamma$ & $\alpha$ & $d_{\eta}$ & $q$ \\\hline
		CPGT-C1 &$C_1$&0.05&0.1&100&0.99\\
		CPGT-C2-1 &$C_2$&0.2&0.1&100&0.99\\
		CPGT-C2-2 &$C_2$&0.05&0.15&100&0.99\\
		DiaDSP &---&---&0.15&100&0.99\\\hline
	\end{tabular}
	\caption{}
		\label{tab:parameter}
\end{table}
\begin{figure}
	\centering
	\includegraphics[width=1\linewidth]{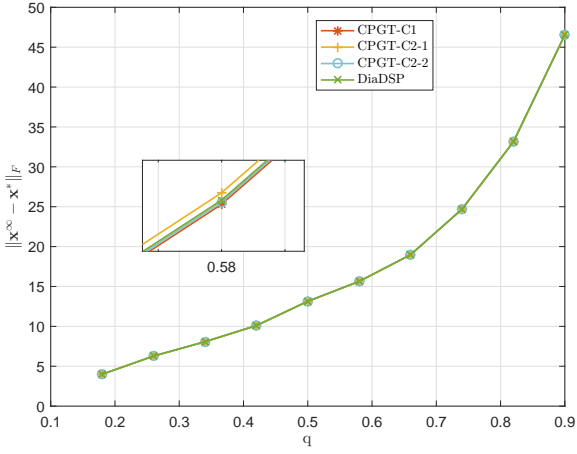}
	\caption{Effect of noise decaying rate on convergence accuracy.}
	\label{fig:accuracy}
\end{figure}

We further simulate the effect of the noise decaying rate on convergence accuracy. Let $\x^*=\mathbf{1}(x^*)^T$ with $x^*\in\R^d$ being the optimal solution. We use $\Vert\x^\infty-\x^*\Vert_F$ to measure the convergence accuracy of different algorithms. We set $d_\eta=5$ and other parameters $\gamma$ and $\alpha$ are the same as Table~\ref{tab:parameter}. The relation between accuracy and decaying rate $q$ is shown in Fig.~\ref{fig:accuracy}, where $q=0.18,0.26,0.34,0.42,0.5,0.58,0.66,0.74,0.82,0.9$. It can be seen that the accuracy of CPGT is nearly the same as that of DiaDSP. Furthermore, accuracy is only noise dependent and not related to stepsize, $\gamma$, and compressors. 
	
\section{conclusion}\label{conclusion}
In this paper, we studied differentially private distributed optimization under limited communication. Specifically, for a class of biased but contractive compressors, we proposed a novel Compressed differentially Private Gradient Tracking algorithm (CPGT). We established linear convergence of CPGT if all the local cost functions are smooth and strongly convex. The proposed CPGT has the same accuracy as the algorithm with idealized communication. Unlike the previous literature, CPGT preserves the differential privacy for the local cost function of each agent with a class of biased but contractive compressors. Future work includes extending to directed graphs and considering more general compressors.


\appendix
\section{The proof of Lemma~\ref{lemma:linearinequalities}}\label{app-linearinequalities}
\subsection{Supporting Lemmas}
\begin{lemma}\label{lemma:stronglyconvex}
Suppose Assumption~\ref{as:smoothandsconvex} holds, if $\alpha\leq\frac{1}{\mu+L}$, one obtains that
\begin{align*}
	\Vert\bx(k)-\x^\infty-\alpha(\bbf(\bx(k))-&\bbf(\x^\infty)\Vert_F\leq\\
	&(1-\alpha \mu)\Vert\bx(k)-\x^\infty\Vert_F.
\end{align*}
\begin{pf}
From the property $\frac{\mathbf{1}\mathbf{1}^T}{n}\bx(k)=\bx(k),~\frac{\mathbf{1}\mathbf{1}^T}{n}\x^\infty=\x^\infty$, we have
\begin{align*}
	\Vert\bx(k)-\x^\infty-&\alpha(\bbf(\bx(k))-\bbf(\x^\infty)\Vert_F\leq\\
	&\Vert\bx(k)-\x^\infty-\alpha(\f(\bx(k))-\f(\x^\infty)\Vert_F.
\end{align*}
Then the proof can be completed by the same line of Lemma 10 in~\cite{qu2017harnessing}.
\end{pf}
\end{lemma}
\begin{lemma}\label{lemma:upperboundoflw}
\cite[]{liao2022compressed} Suppose Assumption~\ref{as:strongconnected} holds. For~$\gamma\in(0,1]$ and any $\omega\in\R^{n\times d}$, we have $\Vert W_\gamma\omega-\bar{\omega}\Vert_F\leq\hat{\lambda}\Vert\omega-\bar{\omega}\Vert_F$, where $\hat{\lambda}=1-\gamma(1-\llw)$ with $\llw=\bar{\lambda}_{W-\frac{\mathbf{1}\mathbf{1}^T}{n}}$.
\end{lemma}

\subsection{The proof of Lemma~\ref{lemma:linearinequalities}}
We prove Lemma~\ref{lemma:linearinequalities} by constructing the upper bounds of $\mathbb{E}[\ox(k+1)],~\mathbb{E}[\obx(k+1)],~\mathbb{E}[\oy(k+1)],~\mathbb{E}[\osx(k+1)]$, and $\mathbb{E}[\osy(k+1)]$, respectively.

(a) According to~\eqref{iterationx1} and Lemma~\ref{lemma:upperboundoflw}, we obtain 
\begin{align}\label{eq:boundofaveandstatex}
	\begin{aligned}
	\E[\ox(k+1)]&\leq\lw\E[\ox(k)]+\alpha\E[\oy(k)]\\
	&\quad+\gamma\bar{\lambda}_{W-I}\E[\osx(k)]+\nu_1(k),
	\end{aligned}
\end{align} 
where $\nu_1(k)=\lw\E[\Vert \mathbf{\eta_x}(k)-\mathbf{\bar{\eta}_x}(k)\Vert_F]$.

(b) From~\eqref{iterationy1}, we first introduce a key property of CPGT, i.e., for $k\geq0$,
	\begin{align}\label{eq:sumofy}
		\mathbf{1}^T\y(k)=\mathbf{1}^T\f(\mathbf{x}(k))+\mathbf{1}^T\sum_{t=0}^{k-1}\ey(t).
	\end{align}
By \eqref{iterationx1},~\eqref{eq:xinfty}, and~\eqref{eq:sumofy}, we have
\begin{align*}
	\bx(k+1)-\x^\infty=&\bx(k)-\x^\infty+\mathbf{\bar{\eta}_x}(k)+\alpha\sum_{t=k}^\infty\mathbf{\bar{\eta}_y}(t)\\
	&\quad-\alpha(\bbf(\x(k))-\bbf(\mathbf{x}^\infty)).
\end{align*}
From Lemma~\ref{lemma:stronglyconvex}, we have
\begin{align}\label{eq:boundofaveandinfty}
	\begin{aligned}
\mathbb{E}[\obx(k+1)]\leq &(1-\alpha\mu)\mathbb{E}[\obx(k)]+\alpha L\E[\ox(k)]+\nu_2(k),
\end{aligned}
\end{align}
where $\nu_2(k)=\mathbb{E}[\Vert\mathbf{\bar{\eta}_x}(k)\Vert_F]+\alpha\sum_{t=k}^\infty\mathbb{E}[\Vert\mathbf{\bar{\eta}_y}(t)\Vert_F]$.

(c) It follows from~\eqref{iterationy1}, Lemma~\ref{lemma:upperboundoflw}, and Assumption~\ref{as:smoothandsconvex} that
\begin{align}\label{eq:boundofy}
	\begin{aligned}
\E[\oy(k+1)]\leq&\lw\E[\oy(k)]+L\E[\Vert\x(k+1)-\x(k)\Vert_F]\\
&+\lw\E[\Vert \mathbf{\eta_y}(k)
-\mathbf{\bar{\eta}_y}(k)\Vert_F]\\
&+\gamma\bar{\lambda}_{W-I}\E[\osy(k)].
\end{aligned}
\end{align}	
By~\eqref{iterationx}, we have
\begin{align}\label{eq:boundofstatex}
\begin{aligned}
	\E[&\Vert\x(k+1)-\x(k)\Vert_F]\leq \gamma\lwi\E[\ox(k)]\\
	&+\gamma\lwi\E[\osx(k)]+\alpha\E[\Vert\y(k)\Vert_F]+\lw\E[\Vert\ex(k)\Vert_F].
\end{aligned}
\end{align}
Furthermore, from~\eqref{eq:xinfty} and~\eqref{eq:sumofy}, one obtains that
\begin{align}\label{eq:normboundofy}
	\begin{aligned}
	\E[\Vert\y(k)\Vert_F]&\leq\E[\oy(k)]+\E[\Vert\by(k)\Vert_F]\\
	&\leq\E[\oy(k)]+\sum_{t=k}^\infty\E[\Vert\eby(t)\Vert_F]+L\E[\ox(k)]\\
	&\quad+L\E[\obx(k)].
	\end{aligned}
\end{align}
Combining~\eqref{eq:boundofy}--\eqref{eq:normboundofy}, it holds that
\begin{align}\label{eq:boundofaveandstatey}
	\begin{aligned}
\E[\oy&(k+1)]\leq L(\gamma\lwi+\alpha L)\E[\ox(k)]\\
&+\alpha L^2\E[\obx(k)]+(\alpha L+\lw)\E[\oy(k)]\\
&+L\gamma\lwi\E[\osx(k)]+\gamma\lwi\E[\osy(k)]+\nu_3(k),
	\end{aligned}
\end{align}
where~$\nu_3(k)=\lw\E[\Vert\ey(k)-\eby(k)\Vert_F]+\alpha L\sum_{t=k}^\infty\mathbb{E}[\Vert\mathbf{\bar{\eta}_y}(t)\Vert_F]$ $+L\lw\E[\Vert\ex(k)\Vert_F]$.

(d) From \eqref{citerationx} and Assumption~\ref{as:compressor}, one obtains that
\begin{align}\label{eq:boundofsigmax}
	\begin{aligned}
\E[\osx(k+1)]&\leq\sqrt{\varphi}\E[\Vert\x(k+1)+\ex(k+1)-\x^c(k)\Vert_F]\\
&\leq\sqrt{\varphi}\E[\osx(k)]+\sqrt{\varphi}\E[\Vert\x(k+1)-\x(k)\Vert_F]\\
&\quad+\sqrt{\varphi}\E[\Vert\ex(k+1)-\ex(k)\Vert_F].
	\end{aligned}
\end{align}
Substituting $\E[\Vert\x(k+1)-\x(k)\Vert_F]$ into \eqref{eq:boundofstatex}--\eqref{eq:normboundofy}, we obtain
\begin{align}\label{eq:boundofsigmaxnorm}
	\begin{aligned}
		\E[\osx(k+1)]&\leq \sqrt{\varphi}(\gamma\lwi+\alpha L)\E[\ox(k)]\\
		&\quad+\sqrt{\varphi}\alpha L\E[\obx(k)]+\sqrt{\varphi}\alpha \E[\oy(k)]\\
		&\quad+\sqrt{\varphi}(\gamma\lwi+1)\E[\osx(k)]+\nu_4(k),
	\end{aligned}
\end{align}
where $\nu_4(k)=\sqrt{\varphi}\lw\E[\Vert\ex(k)\Vert_F]+\!\!\sqrt{\varphi}\E[\Vert\ex(k+1)-\ex(k)\Vert_F]\!+\!\sqrt{\varphi}\alpha\sum_{t=k}^\infty\mathbb{E}[\Vert\mathbf{\bar{\eta}_y}(t)\Vert_F]$.

(e) Similar to~\eqref{eq:boundofsigmax}, we have
\begin{align}\label{eq:boundofsigmay}
	\begin{aligned}
		\E[\osy(k+1)]&\leq\sqrt{\varphi}\E[\osy(k)]+\sqrt{\varphi}\E[\Vert\y(k+1)-\y(k)\Vert_F]\\
		&\quad+\sqrt{\varphi}\E[\Vert\ey(k+1)-\ey(k)\Vert_F].
	\end{aligned}
\end{align}
From \eqref{iterationy}, it holds that
\begin{align}\label{eq:boundofstatey}
	\begin{aligned}
		\E[\Vert\y(k+1)&-\y(k)\Vert_F]\leq \gamma\lwi\E[\oy(k)]\\
		&+\gamma\lwi\E[\osy(k)]+\lw\E[\Vert\ey(k)\Vert_F]\\
		&+L\E[\Vert\x(k+1)-\x(k)\Vert_F].
	\end{aligned}
\end{align}
Combining~\eqref{eq:boundofstatex},~\eqref{eq:normboundofy},~\eqref{eq:boundofsigmay}, and~\eqref{eq:boundofstatey}, we obtain that
\begin{align}\label{eq:boundofsigmaynorm}
	\begin{aligned}
		\E[&\osy(k+1)]\leq \sqrt{\varphi}L(\gamma\lwi+\alpha L)\E[\ox(k)]\\
		&+\sqrt{\varphi}\alpha L^2\E[\obx(k)]+\sqrt{\varphi}(\alpha L+\gamma\lwi) \E[\oy(k)]\\
		&+\sqrt{\varphi}L\gamma \lwi\E[\osx(k)]+\sqrt{\varphi}(\gamma\lwi+1)\E[\osy(k)]\\
		&+\nu_5(k),
	\end{aligned}
\end{align}
where $\nu_5(k)=\sqrt{\varphi}L\lw\E[\Vert\ex(k)\Vert_F]+\!\!\sqrt{\varphi}\E[\Vert\ey(k+1)-\ey(k)\Vert_F]+\sqrt{\varphi}\alpha L \sum_{t=k}^\infty\mathbb{E}[\Vert\mathbf{\bar{\eta}_y}(t)\Vert_F]+\sqrt{\varphi}\lw\E[\Vert\ey(k)\Vert_F]$.

Then the elements of $G$ correspond to the coefficients in \eqref{eq:boundofaveandstatex},~\eqref{eq:boundofaveandinfty},~\eqref{eq:boundofaveandstatey},~\eqref{eq:boundofsigmaxnorm}, and~\eqref{eq:boundofsigmaynorm}. Let $\nu(k)\!\triangleq\![\nu_1(k),\nu_2(k),\nu_3(k),\nu_4(k),\nu_5(k)]^T$\!. Since  $\eta_{x_i}\!\sim\text{Lap}_d(d_{\eta_{x_i}}q_i^k)$ and $\eta_{y_i}\sim\text{Lap}_d(d_{\eta_{y_i}}q_i^k)$, we have $\nu(k)\preceq \vartheta \bar{q}^k\bar{d_\eta}$, where $\vartheta$ is given by
\begin{align}\label{parameters}
\begin{aligned}
	\vartheta=\bigg\{&\lw,~\frac{1-\bar{q}+\alpha}{1-\bar{q}},~(L+1)\lw+\frac{\alpha L}{1-\bar{q}},\\
	&\sqrt{\varphi}(\lw+1)+\frac{\sqrt{\varphi}\alpha }{1-\bar{q}},\\
	&\sqrt{\varphi}L\lw+\sqrt{\varphi}+\frac{\sqrt{\varphi}\alpha L}{1-\bar{q}}+\sqrt{\varphi}\lw \bigg\}n^2.
\end{aligned}
\end{align}

\section{The proof of Lemma~\ref{lemma:convergence}}\label{app-convergence}
\subsection{Supporting Lemma}
We first introduce the following useful lemma.
\begin{lemma}\label{lemma:upperbdofg}
	(Corollary 8.1.29 in~\cite{horn2012matrix}) Let $A\in\R^{n\times n}$ be a nonnegative matrix and let $v\in\R^n$ be an element-wise positive vector. For $\beta\geq0$, $\bar{\lambda}_{G}\leq\beta$ if $Gv\leq\beta v$.
\end{lemma}

\subsection{Proof of Lemma~\ref{lemma:convergence}}
From Lemma~\ref{lemma:upperbdofg}, it is obvious that Lemma~\ref{lemma:convergence} can be proved if the following linear inequalities holds.
\begin{align}\label{eq:boundofg}
	G\zeta\preceq (1-\frac{m}{2}\alpha\mu)\zeta,
\end{align}
where $\zeta=\left[\zeta_1,\zeta_2,L\zeta_3,\zeta_4,L\zeta_5\right]^T$ with $\zeta_1,\dots,\zeta_5>0$.

(a) From~\eqref{eq:boundofaveandstatex}, the first inequality in~\eqref{eq:boundofg} is
\begin{align}\label{eq:boundofg1}
	\lw\zeta_1+\alpha L\zeta_3+\gamma\lwi\zeta_4\leq (1-\frac{m}{2}\alpha\mu)\zeta_1,
\end{align}
which can be rewritten as
\begin{align*}
	\alpha(L\zeta_3+\frac{m\mu}{2}\zeta_1)+\gamma\lwi\zeta_4\leq\gamma(1-\llw)\zeta_1.
\end{align*}
Therefore, inequality~\eqref{eq:boundofg1} holds if $\alpha\leq\frac{\gamma(1-\llw)\zeta_1}{L(2\zeta_3+m\zeta_1)}$ and $\zeta_4\leq\frac{(1-\llw)\zeta_1}{2\lwi}$.

(b) From~\eqref{eq:boundofaveandinfty}, the second inequality in~\eqref{eq:boundofg} is
\begin{align}\label{eq:boundofg2}
\alpha L\zeta_1+(1-\alpha\mu)\zeta_2\leq  (1-\frac{m}{2}\alpha\mu)\zeta_2.
\end{align}
Inequality~\eqref{eq:boundofg2} holds if $\zeta_1\leq\frac{m\mu\zeta_2}{2L}$.

(c) From~\eqref{eq:boundofaveandstatey}, the third inequality in~\eqref{eq:boundofg} is
\begin{align}\label{eq:boundofg3}
\begin{aligned}
L(\gamma\lwi+&\alpha L)\zeta_1+\alpha L^2\zeta_2+(\alpha L+\lw)L\zeta_3\\
&+L\gamma\lwi\zeta_4+L\gamma\lwi\zeta_5\leq  (1-\frac{m}{2}\alpha\mu)L\zeta_3.
\end{aligned}
\end{align}
Diving $\gamma$ and~$L$ on both side of~\eqref{eq:boundofg3}, we have
\begin{align*}
\frac{\alpha }{\gamma}(L\zeta_1+L\zeta_2+L\zeta_3+\frac{m}{2}\mu\zeta_3)&+\lwi(\zeta_1+\zeta_4+\zeta_5)\\
&\leq(1-\llw)\zeta_3.\addtag\label{eq:boundofg31}
\end{align*}
From~\eqref{eq:boundofg31}, inequality~\eqref{eq:boundofg3} holds if $\zeta_1+\zeta_4+\zeta_5\leq\frac{(1-\llw)\zeta_3}{2\lwi}$ and $\alpha\leq\frac{\gamma(1-\llw)\zeta_3}{L(2\zeta_1+2\zeta_2+(2+m)\zeta_3)}$.

(d) From~\eqref{eq:boundofsigmaxnorm}, the fourth inequality in~\eqref{eq:boundofg} is
\begin{align}\label{eq:boundofg4}
	\begin{aligned}
		\sqrt{\varphi}(\gamma\lwi+\alpha L)\zeta_1&+\sqrt{\varphi}\alpha L\zeta_2+\sqrt{\varphi}\alpha L\zeta_3\\
		&+\sqrt{\varphi}(\gamma\lwi+1)\zeta_4\leq  (1-\frac{m}{2}\alpha\mu)\zeta_4.
	\end{aligned}
\end{align}
If~$\alpha\leq\frac{\gamma\lwi}{L}$, one obtains that
\begin{align*}
\begin{aligned}
	\sqrt{\varphi}\gamma\lwi(2\zeta_1+\zeta_2+\zeta_3+\zeta_4)+&\frac{\gamma\lwi m\mu\zeta_4}{2L}\\
	&\leq  (1-\sqrt{\varphi})\zeta_4.
\end{aligned}
\end{align*}
Then, inequality~\eqref{eq:boundofg4} holds if $\alpha\leq\frac{\gamma\lwi}{L}$, $\gamma\leq\min\{\frac{(1-\sqrt{\varphi})\zeta_4}{2\sqrt{\varphi}\lwi(2\zeta_1+\zeta_2+\zeta_3+\zeta_4)},\frac{(1-\sqrt{\varphi})L}{\lwi m\mu}\}$.

(e) From~\eqref{eq:boundofsigmaynorm}, the fifth inequality in~\eqref{eq:boundofg} is
\begin{align}\label{eq:boundofg5}
	\begin{aligned}
		\sqrt{\varphi}L(&\gamma\lwi+\alpha L)\zeta_1+\sqrt{\varphi}\alpha L^2\zeta_2\\
		&+\sqrt{\varphi}(\alpha L+\gamma\lwi)L\zeta_3+\sqrt{\varphi}L\gamma \lwi\zeta_4\\
		&+\sqrt{\varphi}(\gamma\lwi+1)L\zeta_5\leq  (1-\frac{m}{2}\alpha\mu)L\zeta_5.
	\end{aligned}
\end{align}
Similar to~\eqref{eq:boundofg4}, inequality~\eqref{eq:boundofg5} holds if $\alpha\leq\frac{\gamma\lwi}{L}$, $\gamma\leq\min\{\frac{(1-\sqrt{\varphi})\zeta_5}{2\sqrt{\varphi}\lwi(2\zeta_1+\zeta_2+2\zeta_3+\zeta_4+\zeta_5)},\frac{(1-\sqrt{\varphi})L}{\lwi m\mu}\}$. Overall, \eqref{eq:boundofg} holds if $\zeta_1,\dots,\zeta_5$, $\gamma$, and stepsize $\alpha$ satisfy~\eqref{eq:boundofgamma}--\eqref{eq:relationbetwzeta}.

\section{The proof of Theorem~\ref{theo:privacy}}\label{app-privacy}

 From GTPT, it is clear that the observation sequence $\mathcal{H}=\{\mathcal{H}_k\}_{k=0}^\infty$ is uniquely determined by the noise sequences $\ex=\{\ex(k)\}_{k=0}^\infty$, $\ey=\{\ey(k)\}_{k=0}^\infty$, and random sequence~$\mathbf{\varrho}=\{\mathbf{\varrho}(k)\}_{k=0}^\infty$, where $\mathbf{\varrho}(k)\in\R^{n\times d}$ is a matrix and its element $[\mathbf{\varrho}(k)]_{ij}$ is the compression perturbation of $[\sx(k)]_{ij}$. We use function $Z_{\mathcal{F}}$ to denote the relation, i.e., $\mathcal{H}=Z_{\mathcal{F}}(\ex,\ey,\mathbf{\varrho})$, where $\mathcal{F}=\{\x(0),W,\mathcal{S}\}$. From Definition~\ref{def:differentialprivacy}, to show the differential privacy of the cost function $f_{i_0}$, we need to show that the following inequality holds for any observation $\mathcal{H}\subseteq\text{Range}(\text{CPGT})$ and any pair of adjacent cost function sets $\mathcal{S}^{(1)}$ and $\mathcal{S}^{(2)}$,
\begin{align*}
	\mathbb{P}\{(\eta_{x},\eta_{y},&\mathbf{\varrho})\in\Psi|Z_{\mathcal{F}^{(1)}}(\eta_{x},\eta_{y},\mathbf{\varrho})\in\mathcal{H}\}\\
	&\leq e^\epsilon\mathbb{P}\{(\eta_{x},\eta_{y},\mathbf{\varrho})\in\Psi|Z_{\mathcal{F}^{(2)}}(\eta_{x},\eta_{y},\mathbf{\varrho})\in\mathcal{H}\},
\end{align*}
where $\mathcal{F}^{(l)}\!\!=\!\!\{\x(0),W,\mathcal{S}^{(l)}\!\}$, $l\!=\!1,2$, and $\Psi$ denotes the sample space. Then it is indispensable to guarantee $Z_{\mathcal{F}^{(1)}}(\eta_{x},\eta_{y},\mathbf{\varrho})=$ $Z_{\mathcal{F}^{(2)}}(\eta_{x},\eta_{y},\mathbf{\varrho})$, i.e.,
\begin{align*}
	&C(x_i^{a-c,(1)}(k),\varrho(k))=C(x_i^{a-c,(2)}(k),\varrho(k)),\\
	&C(y_i^{a-c,(1)}(k),\varrho(k))=C(y_i^{a-c,(2)}(k),\varrho(k)),
\end{align*}
for $\forall i\in\mathcal{V}$ and any $k\geq0$, where
\begin{align*}
	&x_i^{a-c,(1)}(k)=x_i^{a,(1)}(k)-x_i^{c,(1)}(k-1),\\
	&y_i^{a-c,(1)}(k)=y_i^{a,(1)}(k)-y_i^{c,(1)}(k-1).
\end{align*}
Since $x_i^{c,(l)}(k))=y_i^{c,(l)}(k)=0,k<0,$ $l=1,2,$ from~\eqref{citerationx}--\eqref{citerationy}, we have
\begin{align*}
	x_i^{c,(1)}(k))=x_i^{c,(2)}(k)),
	y_i^{c,(1)}(k))=y_i^{c,(2)}(k)).
\end{align*}
Then one obtains that 
\begin{align}\label{eq:pdfofcompressor}
	\begin{aligned}
&f_c(x_i^{a-c,(1)}(k),\varrho(k))=f_c(x_i^{a-c,(2)}(k),\varrho(k)),\\
&f_c(y_i^{a-c,(1)}(k),\varrho(k))	=f_c(y_i^{a-c,(2)}(k),\varrho(k)),
\end{aligned}
\end{align}
if $x_i^{a,(1)}(k)=x_i^{a,(2)}(k)$ and $y_i^{a,(1)}(k)=y_i^{a,(2)}(k)$, for $\forall i\in\mathcal{V}$. Then due to the property of conditional probability, we have
\begin{align}\label{eq:ineqofprivacy}
\begin{aligned}
&\frac{\mathbb{P}\{(\eta_{x},\eta_{y},\mathbf{\varrho})\in\Psi|Z_{\mathcal{F}^{(1)}}(\eta_{x},\eta_{y},\mathbf{\varrho})\in\mathcal{H}\}}{\mathbb{P}\{(\eta_{x},\eta_{y},\mathbf{\varrho})\in\Psi|Z_{\mathcal{F}^{(2)}}(\eta_{x},\eta_{y},\mathbf{\varrho})\in\mathcal{H}\}}\\
&\leq\frac{\mathbb{P}\{(\eta_{x},\eta_{y},\mathbf{\varrho})\in\Psi|Z_{\mathcal{F}^{(1)}}(\eta_{x},\eta_{y},\mathbf{\varrho})\in\mathcal{H}\}}{\mathbb{P}\{(\eta_{x},\eta_{y},\mathbf{\varrho})\in\Psi|Z_{\mathcal{F}^{(2)}}(\eta_{x},\eta_{y},\mathbf{\varrho})\in\mathcal{H},E_1\}},
\end{aligned}
\end{align}	
where $E_1=\cup_{k=0}^\infty\{x_i^{a,(2)}(k)=x_i^{a,(1)}(k), y_i^{a,(2)}(k)=y_i^{a,(1)}(k),\forall i\in\mathcal{V}\}$ is an event. Thus, from~\eqref{eq:pdfofcompressor}--\eqref{eq:ineqofprivacy}, the proof can be completed in the same way as the proof of Theorem~5 in~\cite{ding2021differentially}.

\bibliography{ref_Antai}        
\end{document}